\magnification=1000
\overfullrule=0pt
\input amssym.def
\input amssym.tex
\normalbaselineskip = 18pt
\normalbaselines
\parskip=8pt plus0pt minus.5pt
\def\p{\,\prec\kern-1.2ex \prec\,}
\def\ch{\raise 0.5ex \hbox{$\chi$}}

\def\nm{{\cal M}}

\def\nmt{\widetilde {\cal M}}

\def\em{E({\cal M},\tau )}
\def\lp{L_p({\cal M},\tau )}
\def\vlp{\Vert _{_{L_p({\cal M},\tau )}}}

\def\lm{L^1({\cal M},\tau )\cap {\cal M}}

\def\mp+{$L^p(\nm)+{\cal M}$}
\def\mp+n{L^p(\nm)+{\cal M}}

\bigskip
\bigskip

\centerline {\bf ON UNBOUNDED $p$-SUMMABLE  FREDHOLM MODULES}
\bigskip
\bigskip
\centerline { {\bf A. L. Carey$^1$, J. Phillips$^2$,  F.A. Sukochev$^{1,3}$}
\footnote{} {Research supported by the Australian Research Council
and NSERC (Canada) \hfil\break
  }}
\footnote{} {\noindent(1)
 School of Mathematical Sciences, University of Adelaide, Adelaide
5005, Australia}
\footnote{} {\noindent 
(2) Department of Mathematics and Statistics,
 University of Victoria, Victoria
BC Canada}
\footnote{}{\noindent (3) Department of Mathematics and Statistics,
Flinders University,
GPO Box 2100, Adelaide, 5001

 E-mails: acarey@maths.adelaide.edu.au, phillips@math.uvic.ca,
sukochev@ist.flinders.edu.au}
\bigskip

\noindent {\bf Abstract } {\it
We prove that odd unbounded $p$-summable Fredholm modules
are also  bounded $p$-summable Fredholm modules (this is the odd
counterpart of a result of A. Connes for the case of even
Fredholm modules). The approach we
use is via
estimates of the form
$$
 \Vert \phi (D)-\phi (D_0)\vlp \leq C\cdot \Vert D-D_0\Vert ^{1\over 2},
$$
where $\phi (t)=t(1+t^2)^{-1/2},\ D_0=D_0^*$ is an unbounded
linear operator affiliated with a semifinite von Neumann algebra ${\cal M}$,
$D-D_0$ is a bounded self-adjoint linear operator from ${\cal M}$ and
$ ({\bf 1}+D_0^2)^{-1/2}\in \lp$, where $\lp$
is a non-commutative $L_p$-space associated with ${\cal M}$. It follows
from our results that, if $p\in (1,\infty )$, then $\phi (D)-\phi (D_0)$
belongs to the space $\lp$.}

\noindent {\bf 0.\ Introduction }

This paper concerns the question arising in the quantised
calculus of Alain Connes [Co1,Co2]
outlined in the abstract. To explain
our results we need some further notation. Let
 ${\cal M}$ be a semifinite von Neumann algebra on a separable Hilbert
space $\cal H$
and let $\lp$ be a non-commutative $L_p$-space associated with
$({\cal M},\tau )$, where $\tau$ is a faithful, normal semifinite trace on
${\cal M}$. The identity in
$\nm $ is denoted by ${\bf 1}$. Let ${\cal A}$ be a unital Banach
$*$-algebra which is represented in ${\cal M}$
via a continuous $*$-homomorphism $\pi$ which,
without loss of generality, we may assume to be faithful.
Where no confusion arises we suppress $\pi$ in the notation.
The fundamental objects of our analysis are explained in the following
definition.

\noindent {\bf Definition 0.1\quad ( [Co1,2], [CP], [S])\quad } {\sl {\bf An
odd unbounded $p$-summable ({\rm respectively}, bounded pre-) Breuer-Fredholm
module for ${\cal A}$}, is a pair $({\cal M},D_0)$
({\rm respectively}, $({\cal M},F_0)$)  where $D_0$( respectively, $F_0$)
is an unbounded ({\rm respectively}, bounded) self-adjoint operator
affiliated with ${\cal M}$ ( {\rm respectively}, in ${\cal M}$) satisfying:
\smallskip
\noindent (1) $({\bf 1}+D_0^2)^{-1/2}$ ({\rm respectively}, 
$|{\bf 1}-F_0^2|^{1/2}$)
 belongs to  $\lp$; and
\hfill\break
\noindent (2) ${\cal A}_0:=\{a\in {\cal A}\ |\ a({\rm dom}D_0)
\subset {\rm dom}D_0, [D_0,a]\in {\cal M}\}$
({\rm respectively},
${\cal A}_p:=\{a\in {\cal A}\ |\ [F_0,a]\in \lp\}$) is a dense
$*$-subalgebra of ${\cal A}$.

\smallskip
\noindent
When $F_0^2={\bf 1}$ we drop the prefix \lq  pre-'.}

\bigskip

\noindent In the special case when $\nm ={\cal B}({\cal H})$ and $\tau $ is the
 standard trace Tr, we shall omit the word \lq\lq Breuer" from the
 definition and speak about {\bf unbounded ({\rm respectively}, bounded)
$p$-summable Fredholm modules}
$({\cal H},D_0)$ ({\rm respectively}, $({\cal H},F_0)$).
In this case, the non-commutative $L_p$-space
coincides with the Schatten-von Neumann ideal $ {\cal C}_p$ of compact
operators and the \lq\lq bounded" part of
Definition 0.1 is a slight extension of Definition 3 from [Co2], p.290 (where
 ${\cal A}_0={\cal A}$ and $F_0^2={\bf 1}$, see also [Co1], Appendix 2
and [Co2], p298). In the
special case when $\nm$ is a semifinite factor, the
\lq\lq unbounded" part of Definition 0.1 coincides with [CP] Definition 2.1;
and, in the case  $\nm ={\cal B}({\cal H})$, it may be considered as
 an odd counterpart of the notion of an unbounded even $p$-summable
pre-Fredholm module from [Co1] (Section 6, Corollary 3 and the remarks
 thereafter).
Definition 0.1 is adapted from [S] where the notion of Fredholm module is
studied in the more general setting of symmetric operator spaces.

The approach of this paper to the study of these
Breuer-Fredholm modules is motivated by
one of the basic problems of perturbation theory
which  may be formulated as follows.

\noindent
{\sl $I$. \quad  If $F$ and $G$ are continuous functions on $(-\infty,\infty )$
under what conditions does the
smallness of $G(D-D_0)$ imply that of $F(D)-F(D_0)$?}

\noindent
We will present a study of  this problem when
the function $F(t)=t(1+t^2)^{-1/2}$, $G(t)=\surd t$
  and   $D_0$ (respectively, $D-D_0$)
is some self-adjoint (respectively, bounded self-adjoint) operator
 affiliated with a semifinite von Neumann algebra ${\cal M}$
(respectively in $\cal M$). We measure the \lq\lq smallness" of
$G(D-D_0)$ (respectively, $F(D)-F(D_0)$)
in the uniform operator norm  (respectively, in the norm of $\lp$).
 This setting
appeared first in [CP] and further extensive considerations were presented
in [S]. The choice of $F$ and $G$ is suggested by the
 theory of unbounded Fredholm modules [Co1,2] and work on  spectral
flow [P1,2], [CP].

\noindent Following Connes' results for the even case (see [Co1] $I$.6),
the difficulties associated with the mapping
$({\cal H},D_0)\longrightarrow ({\cal H},{\rm sgn} (D_0))$ in the odd case
 were outlined in [CP] (see also [S]). Introduce the map $\phi$
defined by
$$
\phi(D)= D({\bf 1}+D^2)^{-1/2}
$$
which is  a smooth approximation of the function sgn 
and hence explains
 our interest in the difference
$
\phi(D)-\phi (D_0).
$
The results presented in this article contribute also
to the study of the mapping
$({\cal M},D_0)\longrightarrow ({\cal M},{\rm sgn} (D_0))$ which was
initiated in [CP] for odd $p$-summable Breuer-Fredholm modules
and continued in [S].
The choice $F=\phi$ is also dictated by the
 following problem suggested from [Co], [CP] and [S].

\noindent
{\sl $II$.\quad  Does it follow from $({\cal M},D_0)$ being an odd unbounded
$p$-summable Breuer-Fredholm module that $({\cal M},{\rm sgn}(D_0))$ is an
odd bounded $p$-summable Breuer-Fredholm module?}

\noindent The major technical problem in our setting for
question $I$ lies in the difference between norms on the right and left
hand sides. In particular, it makes it virtually impossible to apply
well-known double operator integral techniques from [BSo1-3] and
therefore a new technique is required even in the situation
when ${\cal M}$ coincides with the algebra ${\cal B}({\cal H})$ of
all bounded linear operators on ${\cal H}$. Variants
of this new technique are given in [CP] and [S].
In our present approach to problems $I$ and $II$ we
follow the direction outlined in [S], Section 6 where questions
concerning H\"older estimates of the type
$$
 \Vert \phi (D)-\phi (D_0)\vlp \leq C\cdot \Vert D-D_0\Vert^{1/2}\eqno (0.1)
$$
were shown to be relevant to the Lipschitz  continuity of the absolute value
in the setting of non-commutative $L_p$-spaces. The result of [S] Corollary
6.3 (see also Proposition 2.3 in the present paper) asserts that if
$({\bf 1}+D_0^2)^{-1/2}\in \lp$ then the functions
$$
{\Vert D({\bf 1} + D^2)^{-1/2}-
D_{_{0}}({\bf 1} + D_{_{0}}^2)^{-1/ 2}\vlp \over 
\max \{\Vert D-D_0\Vert ^{1/2},
\ \Vert D-D_0\Vert \} }\quad
{\sl and}\quad {\Vert \ (|D|-|D_0|)\cdot
 ({\bf 1} + D_{_{0}}^2)^{-1/ 2}\ \vlp \over 
\max \{ \Vert D-D_0\Vert ^{1/2}, \ \Vert D-D_0\Vert \} }
$$
are bounded or unbounded  simultaneously for each self-adjoint operator
$D-D_0\in {\cal M}$. In other words the question (0.1) is reduced to
the study of the \lq \lq weighted" difference
$$
(|D|-|D_0|)\cdot  ({\bf 1} + D_{_{0}}^2)^{-1/ 2}. \eqno (0.2)
$$
While it is possible to obtain results in the case when $D-D_0$
is $D_0$-bounded (see [CP])
 we will not consider this approach here. 
We use instead
the approach applied in [S] based on the following result obtained
jointly with Yu.B. Farforovskaya.

\bigskip
\noindent {\bf Proposition 0.2\quad (cf. [S] Proposition 6.5)\quad }
{\sl Let $f$ be a Lipschitz function
with constant $1$ and let $p\in [1,2]$. If $T\in {\cal C}_p$ commutes
with $D_0$, then $(f(D)-f(D_0))T\in  {\cal C}_p$
and, moreover}
$$
\Vert (f(D)-f(D_0))T\Vert _{_{{\cal C}_p}}\leq 
\Vert D-D_0\Vert \cdot \Vert T\Vert _{_{{\cal C}_p}}.
$$

\bigskip
\noindent The proof of Proposition 0.2 relies strongly on
the matrix representation of operators from ${\cal C}_p$ and
is not applicable in the case of  $L_p$-spaces affiliated with an arbitrary
 von Neumann algebra.
 In the present article we
shall work with the weighted difference
(0.2) motivated by the methods from [DDPS] where it was established
that the absolute value is Lipschitz continuous in any reflexive $L_p$-space
associated with an arbitrary von Neumann algebra. The following theorem
is our main result. It extends Proposition 0.2 in the special case that $f$
is the absolute value function and
 $p\neq 1$ and contributes further
to the solution of problem $I$.

\noindent {\bf Theorem 0.3\quad } {\sl
(i)  Let $x,y$ be self-adjoint
 operators affiliated with ${\cal M}$ with $x=y+a$ with $a\in {\cal M}$ and
let $z=({\bf 1} + x^2)^{-1/2}\in \lp \cap {\cal M}$ for some fixed
$p\in (1,\infty )$. We then have
$$
(|x|-|y|) z\in  \lp\quad {\sl and} \quad
{y\over ({\bf 1}+y^2)^{1/2}} -  {x\over ({\bf 1}+x^2)^{1/2}}\in \lp.
$$
Moreover
$$
\Vert (|x|-|y|)z\vlp \leq
{\cal Z}_p\max \{ \Vert x-y\Vert ^{1/2}, \ \Vert x-y\Vert \}\cdot \Vert z\vlp \eqno (0.3)
$$
and
$$
\Vert {y\over ({\bf 1}+y^2)^{1/2}} -  {x\over ({\bf 1}+x^2)^{1/2}}\vlp
\leq {\cal Z}_p'\max \{ \Vert x-y\Vert ^{1/2}, \ \Vert x-y\Vert \}
\cdot \Vert z\vlp. \eqno (0.3)'
$$
for some positive constants ${\cal Z}_p$ and ${\cal Z}_p'$ which depend
on $p$ only.

\noindent (ii) Let $x,y$ be  self-adjoint,
$\tau$-measurable operators affiliated with ${\cal M}$.
Let $x=y+a$
with $a\in {\cal M}$ and
suppose that $z\ge 0$ belongs to $\lp \cap {\cal M},$ for some fixed
$p\in (1,\infty )$, commutes with $x$ and has support projection ${\bf 1}$.
Then $(|x|-|y|) z\in  \lp$ and, moreover
$$
\Vert (|x|-|y|)z\vlp \leq
{\cal K}_p\Vert x-y\Vert \cdot \Vert z \vlp
$$
for some positive constant ${\cal K}_p$ which depends on $p$ only.}

\noindent The proofs of theorem 0.3(i) and (ii)
we present here  are independent of each other
(although initially we used (ii) to prove (i)).
Notice the differences in the technical assumptions,
these  are important
and force us to use rather  different arguments.
For the purposes of this paper
 Theorem 0.3(i) is the main result
because from it we can deduce the
following corollary which answers question $II$ in the affirmative
and extends earlier
results of the third named author [S] for the
 case ${\cal M}={\cal B}({\cal H}),\ 1\leq p\leq 2$.

\noindent {\bf Corollary 0.4\quad } {\sl If 
$1<p<\infty$ and  $(\nm , D_0)$ is an odd
 unbounded $p$-summable Breuer-Fredholm module
for the Banach $*$-algebra  ${\cal A}$
 then $(\nm , {\rm sgn}(D_0))$ is an odd  bounded
$p$-summable Breuer-Fredholm module for  ${\cal A}$.}

\noindent
The organisation of the paper is straightforward.
In the next section we shall present a few facts and definitions
 which are necessary for the proof of Theorem 0.3. Our
presentation of the proof of Theorem 0.3 in Section 2 
 requires us to develop further some ideas from 
[DDPS] although  our considerations here are largely
independent of that paper with the exception of one
technical lemma. 
We have deliberately made our discussion
independent of [S] including
the needed results in the Appendix.   Section 3 contains our applications
to non-commutative geometry.

\bigskip
\noindent {\bf 1.\  Preliminaries}\quad We denote
by ${\cal M}$ a semifinite von Neumann algebra on the Hilbert space
${\cal H}$, with a fixed faithful, normal semifinite trace $\tau $.
A linear operator $x$:dom$(x)\to {\cal H} $, with domain
dom$(x)\subseteq {\cal H}$, is said to be {\it affiliated with}  $\nm $  if
$ux=xu$ for all unitaries $u$ in the commutant $\nm'$ of $\nm $
(our basic references for facts about von Neumann algebras are 
[D] or [SZ]).
Given a positive self-adjoint operator $x$ in ${\cal H}$, we denote by
$E^x_t$ (or just $E_t$ 
if there is no danger of confusion)
 the spectral projection of $x$ corresponding to the
interval $(-\infty, t)$. If $x$ is a positive
self-adjoint operator in ${\cal H}$
affiliated with $\nm $, then $E^{ x}_{[0,t)}\equiv E^x_t \in \nm$
 and $xE^{x}_t\in \nm$ for all $t >0$ ([SZ] E.9.10,
E.9.25). If $x$ is a closed linear operator in ${\cal H}$ with polar
decomposition $x=v\vert x\vert$, then $v^*v=s(\vert x\vert)$, where
$s(\vert x\vert)$ is the support projection of $\vert x\vert$ ( [SZ] 9.4).
If $x$ is affiliated with $\nm$, then $v\in \nm$ and $\vert x\vert$ is
affiliated with $\nm$  ([SZ] 9.29). The set of all closed,
densely defined  operators  affiliated with ${\cal M}$
will be denoted by $\widetilde {\widetilde {\cal M}}$.

\noindent An operator $x\in \widetilde {\widetilde {\cal M}}$ is called
$\tau$-{\it measurable} (affiliated with ${\cal M}$) if and only if there
exists $s>0$ such that $\tau (1-E^{\vert x\vert}_s)<\infty$.
The set of all $\tau $-measurable  operators  forms a
$*$-algebra $\widetilde {\cal M}$ with the sum and
product defined as the respective closures of the algebraic sum and
product. For $\epsilon, \delta >0
$ we denote by $N(\epsilon, \delta )$ the set of all $x\in\nmt$ for
which there exists an orthogonal projection $p\in \nm$ such that
$p({\cal H})\subseteq$dom$(x), \Vert xp\Vert \leq \epsilon $ and
$\tau(1-p)\leq\delta$.  The sets $\{
N(\epsilon, \delta ):\epsilon, \delta >0\}$ form a base at $0$ for a
metrizable Hausdorff topology in $\nmt$, which is called the {\it measure
topology}.  Equipped with this measure topology, $\nmt$  is a complete
topological $*$-algebra.  These facts and their proofs can be found in
the papers [Ne], [Te] and [FK]. It is known (see [Ti] and also [DDPS] 
Theorem 1.1) that if $x\in\nmt,\ \{x_n\}_{_{n=1}}^\infty \subset\nmt$ and if
 $x_n\rightarrow x$ for the measure topology, then also 
$|x_n|\rightarrow |x|$ for the measure topology.

\noindent The space $L_{p}({\cal M},\tau ),\ 1\leq p<\infty $
is the Banach space of all operators $A\in \widetilde {\cal M}$ such  that
$\tau (\vert A\vert ^{p})<\infty $ with  the
norm  $\Vert A\vlp:=(\tau (\vert A\vert ^{p}))^{1/p}$,  where
$\vert A\vert =(A^{*}A)^{1/2},\ i=1,2.$ If $\nm ={\cal L}({\cal H})$ and
$\tau $ is the standard trace Tr, then   $\nmt =\nm$ and, then
$L_{_{p}}(\nm ,\tau)$ is precisely the {\it Schatten class}
 ${\cal C}_{_{p}},\ 1\leq p<\infty $.

\noindent
If $\{x_{\alpha }\}_{_{\alpha \in {\cal A}}}\subseteq \nm $ is a net
and if $x\in \nm $, then we will write  
$x _{_{\alpha }}{\buildrel (s)\over
 \rightarrow} x$ to  denote convergence in the $\sigma$-strong 
(operator) topology (see [Ta] p. 68 and [SZ] p.132). If we consider 
the left regular representation of ${\cal M}$ on ${\cal H}=L_2({\cal
M},\tau)$, then it is straightforward that
the convergence in the $\sigma$-strong topology 
coincides with the convergence 
in the strong operator topology. 
It is well-known (see [Da1], p.115, [Si], p. 40 and also [DDPS]
Corollary 1.4)  that if $x _{_{\alpha }}=x _{_{\alpha }}^*,\ \forall \alpha,\ 
\sup_\alpha \Vert x _{_{\alpha }}\Vert _{_{\infty }}<\infty $ and 
if $x _{_{\alpha }} \rightarrow x$ in the strong operator topology, then 
$|x _{_{\alpha }}| \rightarrow |x|$ in the strong operator topology. 
In particular, if $\{e_n\}_{_{n=1}}^\infty$ is a
sequence of projections from ${\cal M}$ such that
$e_n\uparrow _{_{n }}{\bf 1}$ and $x=x^* \in \nm$, then 
$e_n x e_{n } \rightarrow  x$  and $|e_n x e_{n }|\rightarrow  |x|$
 in the strong operator topology, whence
$$
e_n x e_{n }{\buildrel (s)\over \rightarrow}  x\quad {\rm and}\quad
|e_n x e_{n }|{\buildrel (s)\over \rightarrow}  |x|.\eqno (1.1)
$$
In the proof of Theorem 0.3
we shall use the following easily 
verified fact.
If $x _{_{\alpha }}\in \lp$ for $1<p<\infty$ with 
$\Vert x_{_{\alpha }}\vlp \leq C<\infty $ for all $\alpha$ and
either 
 $x_{_{\alpha }} \rightarrow x$ in the measure topology
or we have
$x_{_{\alpha }}=x _{_{\alpha }}^*$ for all $\alpha$,  
$\sup_\alpha \Vert x _{_{\alpha }}\Vert _{_{\infty }}<\infty$  and 
$x_{_{\alpha }}{\buildrel (s)\over \rightarrow} x,$  then 
$$
x\in \lp\quad {\rm and}\quad \Vert x\vlp \leq C.\eqno (1.2)
$$ 
The rigorous proof of the latter fact in a slightly more 
general situation may be found in [DDPS] Proposition 1.6 and in [FK] 
Theorems 3.5, 3.6.

\bigskip
\noindent An important fact from the geometry of non-commutative $L_p$-spaces
used in [DDPS] is that any reflexive $L_p$-space associated with an
arbitrary semifinite von Neumann algebra $({\cal M},\tau )$ is a UMD-space
(see [BGM1]). An equivalent form of the latter fact is that the $\lp$-valued
 generalization of the Riesz projection is bounded in any
Bochner space $L_p(G, \lp)$, where $G$ is an arbitrary connected compact
Abelian group, the Riesz projection is defined with respect to a positive
cone of a linear ordering of the dual group $\hat G$ and $p\in (1,\infty )$.
This fact together with the so-called transference method (see [BGM2]) was
used in [DDPS] to establish the following result (which in the special case 
${\cal M}={\cal B}({\cal H})$ 
was first established by E.B. Davies in [Da2]).

\noindent {\bf Lemma 1.1\quad (cf. [DDPS] Lemma 3.2) }\quad {\sl
If $1<p<\infty $, then there
exists a constant $K_p>0$, which depends only on $p$ such that}
$$
\left \Vert \sum _{m,n=1}^N{\lambda _m-\mu _n\over \lambda _m+\mu
_n}p_map_n
\right \vlp \leq K_p\Vert a\vlp,
$$
{\sl for all semifinite von Neumann algebras $(\nm ,\tau )$, for all
finite sequences $p_1,p_2, \dots ,p_N$ of mutually orthogonal
projections in $\nm $, for all $a\in
\lp $ and all choices $0\leq \lambda _1,\lambda _2,\dots
,\lambda _N;\ \mu _1,\mu _2, \dots ,\mu _N\in {\Bbb R}$ with $\lambda _m
+\mu _n>0$ for all} $m,n=1,2, \dots ,N$.

\bigskip
\noindent {\bf 2.\  Lipschitz and commutator estimates}\quad
This Section contains the main proofs. We begin with
three technical Propositions. The first gives an
estimate of the \lq\lq weighted" commutator $[x,y]z$ which generalizes
similar considerations of [DDPS].

\noindent {\bf Proposition 2.1 }\quad {\sl If $x=x^*\in \lm$, if
$y\in {\cal M}$
and if $z\in \lp\cap \nm$ commutes with $x$, then }
$$
\Vert \ [\vert x\vert ,y]z\ \vlp \leq 2(1+K_p)\Vert \ [x,y]z \ \vlp .
$$

\noindent {\bf Proof }\quad
 Let $x\in \nm \cap L_1(\nm ,\tau )$
be a self-adjoint element of the form
$$
x=(\lambda _1 p_1 +\lambda _2 p_2 +\dots +\lambda _N p_N)-
(\mu _1q_1+\mu _2q_2 +\dots +\mu _N q_N)
$$
where $p_1,p_2,\dots ,p_N,q_1,q_2, \dots ,q_N$ are mutually
orthogonal projections in $\nm $ and $\{\lambda _i\}_{_{i=1}}^N,\
\{\mu _j\}_{_{j=1}}^N\subset [0,\infty )$.
Note that there is no loss of generality in having the same number
of $p_i$'s and $q_j$'s as we can allow some of them to be zero.
Let $z\in \lp\cap \nm$ commute
with these projections.
It follows immediately from Lemma 1.1 that
$$
\left \Vert \sum _{m,n=1}^N{\lambda _m-\mu _n\over \lambda _m+\mu
_n}p_myp_nz
\right \Vert _{_{p}}\leq K_p\Vert yz\Vert _{_{p}},\eqno (2.1)
$$
Letting
$$
y':=(\sum _{m=1}^N\lambda _mp_m)y+y(\sum _{n=1}^N\mu _np_n),\quad y''=(\sum _{m=1}^Np_m)y'(\sum _{n=1}^Np_n)
$$
we see that
$$
p_my''p_n=p_my'p_n=(\lambda _m+\mu _n)p_myp_n,\quad m,n=1,2, \dots ,N
$$
and we infer from (2.1) applied to $y''$ and $z$ that
$$
\eqalignno {
\Vert \sum _{m,n=1}^N( \lambda _m-\mu_n)p_myp_nz \Vert _{_{p}} & =
\left \Vert \sum _{m,n=1}^N{\lambda _m-\mu _n\over \lambda _m+\mu
_n}p_my''p_nz\right \Vert _{_{p}}\cr &
\leq  K_p\Vert y''z\Vert _{_{p}}\cr &
=\Vert \sum _{m,n=1}^N( \lambda _m+\mu_n)p_myp_nz \Vert _{_{p}}. &{(2.2)}\cr  }
$$

 We set
$$
p=\sum _{i=1}^Np_i,\quad q=\sum _{j=1}^Nq_j
$$
and note that without loss of generality we may assume that $p+q={\bf 1}$.
Following [DDPS] Proposition 2.4 (vii) $\Longrightarrow$ (viii))
we have, as $z$ commutes with $p$ and $q$,
 $$
\eqalign {
p[\vert x \vert ,y]zp=\sum _{i=1}^N\sum _{j=1}^N(\lambda
 _i-\lambda _j)p_iyp_jz=&p[x,y]zp\cr
 p[\vert x \vert ,y]zq=\sum _{i=1}^N\sum _{j=1}^N(\lambda _i-\mu
 _j)p_iyq_jz,\quad& p[x,y]zq=\sum _{i=1}^N\sum _{j=1}^N(\lambda _i+\mu
 _j)p_iyq_jz,\cr
 q[\vert x \vert ,y]zp=\sum _{i=1}^N\sum _{j=1}^N(\mu _j-\lambda
 _i)q_jyp_iz,\quad& q[x,y]zp=-\sum _{i=1}^N\sum _{j=1}^N(\mu _j+\lambda
 _i)q_jyp_iz,\cr
 q[\vert x \vert ,y]zq=\sum _{i=1}^N\sum _{j=1}^N(\mu _i-\mu _j)q_iyq_jz
=&- q[x,y]zq. \cr }
$$
Using (2.2) we now have
$$
\eqalign {
\Vert p[\vert x\vert ,y]zp\vlp =\Vert p[x,y]zp\vlp ,\quad
&\Vert p[\vert x\vert ,y]zq\vlp \leq K_p\Vert p[x,y]zq\vlp ,\quad \cr
\Vert q[\vert x\vert ,y]zp\vlp \leq K_p\Vert q[x,y]p\vlp ,\quad
&\Vert q[\vert x\vert ,y]zq\vlp =\Vert q[x,y]zq\vlp ,\quad \cr }
$$
It now follows that
$$
\Vert \ [\vert x\vert ,y]z\ \vlp =\Vert (p+q)[\vert x\vert ,y]z(p+q)\vlp
\leq 2(1+K_p)\Vert\ [x,y]z\ \vlp.\eqno (2.3)
$$

\noindent We suppose now that
$x$ is an arbitrary self-adjoint element from $ \nm \cap L_1(\nm ,\tau )$.
There exists a sequence $\{x_n\}\in \lm $ such that each $x_n,n\geq
1$ is a finite linear combination of spectral projections of $x$, such that
$x_n\to x,\vert x_n\vert \to \vert x\vert $ in $\lm $. It follows
that $[x_n,y]\to [x,y],\ [\vert x_n\vert ,y]\to [\vert x\vert ,y]$ in
$\lm $ and hence
$[x_n,y]z\to [x,y]z,\ [\vert x_n\vert ,y]z\to [\vert x\vert ,y]z$ in
$\lp $ by the continuity of the embedding of
$\lm $ into $\lp $ (here we have adopted the argument used in [DDPS]
Proposition 2.4 (viii) $\Longrightarrow$ (ii)). Noting that $x_n$ commutes with $z$ for every $n\ge 1$ we have via (2.3) that
$$
\Vert \ [\vert x_n\vert ,y]z\ \vlp \leq 2(1+K_p)\Vert \ [x_n,y]z\ \vlp
$$
for all $n\geq 1$ (note that the assumptions
$x_n\in \lm,\ n=1,2,\dots $ and $y\in {\cal M}$ guarantee $[x_n,y]z\in \lp$),
and so also
$$
\Vert \ [\vert x\vert ,y]z\ \vlp \leq 2(1+K_p)\Vert \ [x,y]z \ \vlp .
$$
This completes the proof of Proposition 2.1.\quad $\square$

\noindent We shall now modify a matrix argument from the proof
of the implication (ii)$\Rightarrow $(i) in [DDPS] Theorem 2.2.
It should be noted that the assumptions imposed on the
 element $y$ in the next Proposition are
more stringent than in Proposition 2.1.

\noindent {\bf Proposition 2.2 }\quad {\sl If $x=x^*,y=y^*\in \lm$  and
$z\in \lp\cap \nm$ commutes with $x$, then }
$$
\Vert (\vert x\vert -\vert y\vert)z \vlp \leq 2(1+K_p)\Vert (x-y)z\vlp .
$$

\noindent {\bf Proof }\quad  It should be noted now that the assertion of
Proposition 2.1 holds for
an arbitrary semifinite von Neumann
algebra, in particular it holds for the von Neumann
algebra $\nm _1:=\nm \otimes M_2({\Bbb C})$ of all $2\times 2$ matrices
$$
[x_{ij}]=\pmatrix {x_{11}&x_{12}\cr x_{21}&x_{22}\cr }
$$
with $x_{ij}\in \nm , \ i,j=1,2$, acting on the Hilbert space
${\cal H}\oplus {\cal H}$ with the trace $\tau _1$ given by setting
$$
\tau _1([x_{ij}])=\tau (x_{11})+\tau (x_{22}).
$$
If
$$
X^*=X=\pmatrix {x&0\cr 0&y\cr },\quad Y=\pmatrix {0&0\cr 1&0\cr },\quad
Z=\pmatrix {z&0\cr 0&0\cr },
$$
then
$$
\left [\left \vert
X\right \vert ,Y\right ]Z=
\pmatrix {0&0 \cr (\vert y\vert -\vert x\vert)z &0\cr },
$$
and
$$
\left [X,Y\right ]Z=
\pmatrix {0&0 \cr ( y - x)z &0\cr }.       \eqno (2.4)
$$
Since $X=X^*\in L_1(\nm _1,\tau _1)\cap \nm _1,\ Y\in {\cal M}_1$ and
$ZX=XZ,\ Z\in L_p({\cal M}_1, \tau _1)\cap {\cal M}_1$,
it follows from Proposition 2.1 and (2.4) that
$$
\Vert \ [\left \vert
X\right \vert ,Y ]Z\ \Vert _{_{L_p({\cal M}_1,\tau_1 )}}
 \leq 2(1+K_p)\Vert \ [X ,Y ]Z\ \Vert _{_{L_p({\cal M}_1,\tau_1 )}}
$$
{}From (2.4) it is clear that
$$\Vert \ [\left \vert
X\right \vert ,Y ]Z\ \Vert _{_{L_p({\cal M}_1,\tau_1 )}}
= \Vert (\vert y\vert -\vert x\vert)z\vlp, \quad
\Vert \ [X ,Y ]Z\ \Vert _{_{L_p({\cal M}_1,\tau_1 )}}
= \Vert ( y - x)z\vlp
$$
and the assertion of the proposition follows. $\square$

\noindent Our proof of Theorem 0.3(i) rests on
Proposition 2.3 below combined with a refinement of the approach from [S].
Crucial to our arguments is
the following inequality (Theorem 6.2 of [S])
whose proof we include in the Appendix so that this
paper may be read independently of [S]. For any
$x=x^*,y=y^*\in \widetilde {\widetilde {\cal M}}$ such that
$ x=y+a$ with $a \in \nm$ and  $z=({\bf 1} + x^2)^{-1/2}\in \lp$ we have
$$
 \Vert {|y|\over (1+y^2)^{1/2}}-{|x|\over (1+x^2)^{1/2}}
 \vlp \leq 2^{3/2} \Vert z \vlp \cdot \max \{\Vert x-y\Vert ^{1/2},
\ \Vert x-y\Vert \}.
\eqno (2.5)
$$
Our final technical
 Proposition rests on (2.5) and  is a slight refinement of
[S] Corollary 6.3.

\noindent {\bf Proposition 2.3 }\quad {\sl Let
$x=x^*,y=y^*\in \widetilde {\widetilde {\cal M}}$ and  let
$z=({\bf 1} + x^2)^{-1/2}\in \lp$. If $x=y+a$ with $a\in \nm$
 and the following inequality
holds for some constant $c_p>0$
$$
\Vert \ (|y|-|x|)z\ \vlp \leq
c_p \max \{ \Vert x-y\Vert ^{1/2}, \ \Vert x-y\Vert \}  \eqno (2.6)
$$
then
we have
$$
\Vert y({\bf 1} + y^2)^{-1/2}-
x({\bf 1} + x^2)^{-1/ 2}\vlp \leq c'_p\max \{\Vert x-y\Vert ^{1/2},
\ \Vert x-y\Vert \}. \eqno (2.7)
$$
for the constant
$c'_p:=c_p+  (2^{3/2}+1) \Vert z \vlp$

In its turn, if (2.7) holds for some
constant $c'_p$ and some (self-adjoint) $a=x-y\in \nm$, then (2.6) holds
 with $c_p:=c'_p+  (2^{3/2}+1) \Vert z \vlp$.}

\noindent {\bf Proof \quad} Let (2.6) be satisfied for some
constant $c_p$ and
 some $a=x-y\in \nm$. Then from the equality
$$
|y|\bigl ({{\bf 1}\over ({\bf 1}+y^2)^{1/2}}- {{\bf 1}\over ({\bf 1}+x^2)^{1/2}}\bigr )=
-(|y|-|x|){{\bf 1}\over ({\bf 1}+x^2)^{1/2}} +
 {|y|\over ({\bf 1}+y^2)^{1/2}}-{|x|\over ({\bf 1}+x^2)^{1/2}}
$$
and (2.5) we infer that
$$
\eqalign {
\Vert \ |y| \bigl ({{\bf 1}\over ({\bf 1}+y^2)^{1/2}}- {{\bf 1}\over ({\bf 1}+x^2)^{1/2}}\bigr )\vlp
& \leq \Vert (|y|-|x|)z\vlp +
\Vert {|y|\over ({\bf 1}+y^2)^{1/2}}-{|x|\over ({\bf 1}+x^2)^{1/2}} \vlp \cr &
\leq (c_p +
2^{3/2} \Vert z \vlp)\max \{ \Vert x-y\Vert ^{1/2}, \ \Vert x-y\Vert \}.\cr }
$$
It follows immediately that
$$
\Vert  y \bigl ({{\bf 1}\over ({\bf 1}+y^2)^{1/2}}- {{\bf 1}\over ({\bf 1}+x^2)^{1/2}}\bigr )\vlp \leq
(c_p+  2^{3/2} \Vert z \vlp)\max \{ \Vert x-y\Vert ^{1/2}, \ \Vert x-y\Vert \}.\eqno (2.8)
$$
Now from the equality
$$
y\bigl ({{\bf 1}\over ({\bf 1}+y^2)^{1/2}}-{{\bf 1}\over ({\bf 1}+x^2)^{1/2}}\bigr )=
{y\over ({\bf 1}+y^2)^{1/2}}-{x\over ({\bf 1}+x^2)^{1/2}}+
(y-x){1\over (1+x^2)^{1/2}}
$$
combined with (2.8) we arrive at (2.7)
$$
\eqalign {
\Vert {y\over ({\bf 1}+y^2)^{1/2}}-{x\over ({\bf 1}+x^2)^{1/2}}\vlp & \leq
(c_p+  2^{3/2} \Vert z \vlp)\max \{ \Vert x-y\Vert ^{1/2}, \ \Vert x-y\Vert \} + \Vert z \vlp \cdot \Vert x-y\Vert \cr &
\leq (c_p+  (2^{3/2}+1) \Vert z \vlp)\max \{ \Vert x-y\Vert ^{1/2}, \ \Vert x-y\Vert \}.\cr }
$$
The second assertion is established similarly.\quad $\square$

\noindent With these preliminary results
established we are now in a position to prove Theorem 0.3.

\noindent {\bf Proof of Theorem 0.3(i)\quad }
  Recall that we assume that
$x=x^*,y=y^*\in \widetilde {\widetilde {\cal M}}$ are such that
$ x=y+a$ with $a \in \nm$ and that
$z=({\bf 1} + x^2)^{-1/2}\in \lp$. Introduce the sequence 
$\{e_n=E^z_{[1/n,1]}\}_{n=1}^\infty \subset {\cal M}$.
Note that
$$
e_nz=ze_n,\quad \forall n\ge 1,\quad e_n\uparrow _{_{n }}{\bf 1}.
$$
It is straightforward to see
 from the definition of $e_n,\ n\ge 1$
and the fact $z=({\bf 1} + x^2)^{-1/2}$, that
$$
e_nx=xe_n \in \nm,\quad \forall n\ge 1
$$
and that $e_n\leq nz\in\lp$.
It is immediate that $e_n\in{L_1({\cal M},\tau )}$ for all $n\ge 1$
 and further that 
$$e_nxe_n,\ e_nye_n \in{L_1({\cal M},\tau)\cap {\cal M}}$$
 for all $n\ge 1$.


\noindent  Appealing to proposition 2.2 we have
$$
\eqalignno {
\Vert (\vert e_nxe_n\vert -\vert e_nye_n\vert)e_nze_n \vlp
& \leq {\cal K}_p\Vert e_n(x-y)e_nze_n\vlp \cr &
\leq {\cal K}_p \Vert z\vlp \max \{ \Vert e_n(x-y)e_n\Vert ^{1/2}, \ 
\Vert e_n(x-y)e_n\Vert \},&{(2.9)}\cr }
$$
for all $ n\ge 1$ and all (self-adjoint) $x-y\in \nm$. Noting that
$$
\eqalign {
(\vert e_nxe_n\vert -\vert e_nye_n\vert)e_nze_n & = (\vert e_nxe_n\vert -
\vert e_nye_n\vert)e_n ({\bf 1}+x^2)^{-1/2}e_n\cr &
= (\vert e_nxe_n\vert -\vert e_nye_n\vert)(e_n+(e_nxe_n)^2)^{-1/2}\cr &
=(\vert e_nxe_n\vert -\vert e_nye_n\vert)e_n({\bf 1}+(e_nxe_n)^2)^{-1/2}\cr &
=(\vert e_nxe_n\vert -\vert e_nye_n\vert)({\bf 1}+(e_nxe_n)^2)^{-1/2} \cr }
$$
we may now combine (2.9) with Proposition 2.3 to obtain
$$
\eqalignno {
\Vert {e_nye_n\over ({\bf 1}+(e_nye_n)^2)^{1/2}} -  {e_nxe_n\over ({\bf 1}+(e_nxe_n)^2)^{1/2}}\vlp &
\leq c_p'\max \{ \Vert e_nxe_n-e_nye_n\Vert ^{1/2}, \ \Vert e_nxe_n-e_nye_n\Vert \} \cr &
\leq c_p'\max \{ \Vert x-y\Vert ^{1/2}, \ \Vert x-y\Vert \} &{(2.10)}\cr }
$$
for all $n\ge 1$ and all (self-adjoint) $x-y\in \nm$ with
$$
c_p'=({\cal K}_p+2^{3/2}+1)\Vert z\vlp.\eqno (2.11)
$$

\noindent It should be noted that since $e_n\uparrow _{_{n }}{\bf 1}$ and
$y=x-a, \ a \in \nm$ we have
$$
e_nxe_n(\xi) \to x(\xi),\quad e_nye_n(\xi) \to y(\xi),
$$
as $n\to \infty$ for any $\xi \in$ dom$x$=dom$y$.

\noindent Combining [RS1] Theorems $VIII$.25(a) and $VIII$.20(b) we have
$$
{e_nye_n\over ({\bf 1}+(e_nye_n)^2)^{1/2}} -
{e_nxe_n\over ({\bf 1}+(e_nxe_n)^2)^{1/2}} \
{\longrightarrow}\  {y\over ({\bf 1}+y^2)^{1/2}} -
{x\over ({\bf 1}+x^2)^{1/2}},\quad n\to \infty \eqno (2.12)
$$
in the strong operator topology as $n\to \infty$.

\noindent  We may now proceed to the final part of the
proof. All the operators in 
(2.12) are  uniformly bounded. Hence we may 
apply (1.2) to deduce that,
$$
\Vert {y\over ({\bf 1}+y^2)^{1/2}} -
 {x\over ({\bf 1}+x^2)^{1/2}}
\vlp\leq \liminf_{_{n}} \Vert {e_{n}ye_{n}\over
({\bf 1}+(e_{n}ye_{n})^2)^{1/2}} -
{e_{n}xe_{n}\over
({\bf 1}+(e_{n}xe_{n})^2)^{1/2}}\vlp
$$
and, using (2.10), 
$$
\Vert {y\over ({\bf 1}+y^2)^{1/2}} -
 {x\over ({\bf 1}+x^2)^{1/2}}
\vlp
\leq c_p'\max \{ \Vert x-y\Vert ^{1/2}, \ \Vert x-y\Vert \}.
$$
Letting  (see equality (2.11))
$$
{\cal Z}_p':={\cal K}_p+2^{3/2}+1
$$
we arrive at the inequality $(0.3)'$.
The inequality (0.3) of Theorem 0.3(i) follows from the inequality $(0.3)'$
 via Proposition 2.3.\quad $\square$

\bigskip
\noindent {\bf Proof of Theorem 0.3(ii)\quad }
We suppose first that $x=x^*,y=y^*\in \nm$. Let
$$
e_n:={\bf 1}- E^z_{1/n},\quad  n=1,2,\dots\ .
$$
Then, using $0\leq z\in \lp\cap \nm$ and $s(z)={\bf 1}$, we have that 
$e_n\uparrow _{_{n }}{\bf 1}$ and that
$$
e_nz=ze_n,\quad e_nx=xe_n,\quad
\tau (e_n)<\infty, \quad \forall n\ge 1.
$$
Since $e_n\uparrow _{_{n }}{\bf 1}$ we have (see (1.1))
$$
e_n x e_{n }{\buildrel (s)\over \rightarrow} x ,
\quad e_n y e_{n }{\buildrel (s)\over \rightarrow} y ,
\quad   \vert e_n x e_{n }
\vert {\buildrel (s)\over \rightarrow}\vert x \vert ,
\quad   \vert e_n y e_{n }
\vert {\buildrel (s)\over \rightarrow} \vert y \vert. \eqno (2.13)
$$
It follows from the assumptions $x,y\in \nm $ and from the inequality
$\tau (e_n)<\infty$ that
$$
e_nxe_n=(e_nxe_n)^*,\
e_nye_n=(e_nye_n)^* \in \lm, \quad \forall n\ge 1.
$$
Since
$$
e_n(x-y)ze_n=e_n(x-y)e_nz\quad {\rm and}\quad
\Vert e_n(x-y)ze_n\vlp \leq \Vert (x-y)z\vlp,
$$
it follows from Proposition 2.2 that
$$
\Vert (\vert e_nxe_n\vert -\vert e_nye_n\vert)z \vlp \leq 2(1+K_p)\Vert e_n(x-y)e_nz\vlp
\leq 2(1+K_p)\Vert (x-y)z\vlp, \quad \forall n\ge 1.\eqno (2.14)
$$
The inequality
$$
\Vert (\vert x\vert -\vert y\vert)z \vlp
\leq 2(1+K_p)\Vert (x-y)z\vlp.\eqno (2.15)
$$
is now clear from (2.13), (2.14) combined with (1.2).

\noindent We shall assume now that $x=x^*,y=y^*\in \nmt,\ x-y \in \nm$.

\noindent There exist self-adjoint
projections $\{p_n\}\subseteq \nm $ such that $p_n\uparrow _n {\bf 1},
\tau ({\bf 1}-p_n)\to 0$ and such that $xp_n\in \nm$, $ n\geq 1$. It follows
immediately from the assumption $x-y \in \nm$ that  $yp_n\in \nm, n\geq 1$.
Since it is possible to choose the sequence $\{p_n\}_{n=1}^\infty$
from the set of spectral projections of $x$ we also have
$$
p_{n }z=ze_{n },\quad p_{n }x=xp_{n },\quad n\geq 1.
$$
Thus, appealing to the preceding part of the proof and applying (2.15) we get
$$
\eqalignno {
\Vert (\vert e_{n }xe_{n }\vert -\vert e_{n }ye_{n }\vert)z \vlp
& \leq 2(1+K_p)\Vert (e_{n }xe_{n }-e_{n }ye_{n })z\vlp \cr & =
2(1+K_p)\Vert e_{n }(x-y)ze_{n }\vlp \cr &
\leq 2(1+K_p)\Vert (x-y)z\vlp. &{(2.16)}\cr }
$$
It is easily seen that $p_nxp_n\to x, p_nyp_n\to y$ for the
measure topology and therefore (see Section 1) we have
$$
\vert p_nxp_n\vert -\vert p_n y p_n \vert \to
\vert x \vert -\vert y\vert
$$
for the measure topology. It follows immediately that
$$
(\vert p_nxp_n\vert -\vert p_n y p_n \vert)z \to
(\vert x \vert -\vert y\vert)z
$$
for the measure topology. This fact, combined with (2.16) and (1.2) implies
$$
\Vert (\vert x\vert -\vert y\vert)z \vlp
\leq 2(1+K_p)\Vert (x-y)z\vlp.
$$
The proof of Theorem 0.3(ii) is completed with ${\cal K}_p=2(1+K_p). \quad \square$

\noindent{\bf 3. Applications}

Before we move to the applications we need a preliminary result.
Recall that $\cal A$ is
 a Banach$*$-algebra with a bounded $*$-representation
$\pi:{\cal A}\rightarrow {\cal M}$.
 As the kernel of this $*$-representation is a closed two-sided 
$*$-ideal in {\cal A}, we can (and do) assume for our purposes that $\pi$ is 
faithful. We let $||.||_{\cal A}$ denote the Banach$*$-algebra norm on $\cal A$
 and by renorming $\cal A$ if necessary we can (and do) assume
that 
$ ||\pi(a)|| \leq||a||_{\cal A}$. From now on we
suppress the notation $\pi$, but not the distinct norm $||.||_{\cal A}$ on 
$\cal A$.

\bigskip
\noindent{\bf Lemma 3.1} {\it
The set $${\cal A}_0:=\{a\in {\cal A}\ |\ a({\rm dom}D_0)
\subset {\rm dom}D_0, [D_0,a]\in \nm\}$$
is a Banach $*$-algebra  in the
norm $||a||_0= ||a||_{\cal A}
+||[D_0,a]||$.}
\bigskip

\noindent{\bf Proof}. This result appears to be well known ({\it cf} [CM])
and in any case is a good exercise in careful applications of the definition
of the adjoint of an unbounded operator, see [RN] pages 299-300. $\square$

\noindent We are now in a position to present the proof of Corollary 0.4.
We follow the argument in  [S] Corollary 6.8.

\noindent {\bf Proof of Corollary 0.4}

\noindent We first show that
$(\nm,\phi (D_0))$ is an odd  bounded 
$p$-summable pre-Breuer-Fredholm module for ${\cal A}$.
 Recall that
by assumption (see the part (1) of Definition 0.1 for unbounded
 odd Breuer-Fredholm modules) the element
 $({{\bf 1}\over {\bf 1}+D_0^2})^{1/2}$ belongs to $\lp$ and therefore
$$
({\bf 1}-\phi (D_0)^2)^{1/2}=
({\bf 1}-{D_0^2\over {\bf 1}+D_0^2})^{1/2}=({{\bf 1}\over {\bf 1}+D_0^2})^{1/2}
$$
belongs to $\lp$ too. Thus part (1)
 of Definition 0.1 for bounded odd pre-Breuer-Fredholm modules is satisfied. Thus,
we need to check only the second part of Definition 0.1. It suffices to
show that
$$
{\cal A}_0\subseteq {\cal A}_p.\eqno (3.1)
$$
Using lemma 3.1 we may now apply a result
of [Pa], Theorem 7, to see that the linear span of the set of all unitary
elements $U({\cal A}_0)$ coincides with ${\cal A}_0$.  Hence
in order to establish (3.1) we need to show only
that
$$
[\phi (D_0),u]\in \lp,\quad \forall u\in U({\cal A}_0).\eqno (3.2)
$$
To establish (3.2), we note that for an arbitrary $u\in U({\cal A}_0)$, we have
$$
[\phi (D_0),u]=\phi (D_0)u-u\phi (D_0)=u(u^*\phi (D_0)u-\phi (D_0))=
u(\phi (u^*D_0u)-\phi (D_0)).
$$
{}From our assumptions (see the part (2) 
of Definition 0.1 for unbounded odd Breuer-Fredholm modules) we have that
$$
u^*D_0u-D_0=u^*[D_0,u] \in \nm,\quad \forall u\in U({\cal A}_0).
$$
Therefore, letting $D=u^*D_0u$, we have by Theorem 0.3(i), that
$$
\phi (D)-\phi (D_0)\in \lp
$$
and this shows immediately that (3.2) holds.

\noindent It is now easy to verify that $(\nm,{\rm sgn}(D_0))$ is an odd
 bounded $p$-summable Breuer-Fredholm module for ${\cal A}$.
Indeed,  condition (1) from Definition 0.1 obviously holds. To verify
that condition (2) holds, we note that
$$
\eqalign {
({\rm sgn}(D_0)-\phi (D_0))({\rm sgn}(D_0)+\phi (D_0))
&={\bf 1}-\phi (D_0)^2 \cr
&={\bf 1}-D_0^2(1+D_0^2)^{-1}\cr
&=({\bf 1}+D_0^2)^{-1}\cr
&\leq ({\bf 1}+D_0^2)^{-1/2}\in \lp,\cr}
$$
and since
$$
\bigl ({\rm sgn}(D_0)+\phi (D_0) \bigr)^{-1}\in \nm
$$
it follows that
$$
{\rm sgn}(D_0)-\phi (D_0)=
({\bf 1}+D_0^2)^{-1}\bigl (({\rm sgn}(D_0)+\phi (D_0)\bigr)^{-1}\in \lp,
$$
whence (via the first part of the proof)
$$
[{\rm sgn}(D_0),u]=[{\rm sgn}(D_0)-\phi (D_0),u]+[\phi (D_0),u]\in \lp
$$
for any $u\in U({\cal A}_0)$.\quad $\square$

\noindent The significance of this result for Connes' quantised
calculus is that it fills a lacuna in [Co1]. There the relationship
between bounded and unbounded Fredholm modules is presented
in the even case but not in the odd case. This is
rectified by taking a different viewpoint  in [Co2]
utilising the Dixmier trace, a device which is clearly natural from
the viewpoint of the geometric examples described there.
However using Corollary 0.4 we can fill this lacuna
by a different method which we now explain.

Given an odd $p$-summable unbounded Breuer-Fredholm module
$({\cal M}, D_0)$ for
the algebra $\cal A$,  we have
by Corollary 0.4, that $({\cal M}, F_0={\rm sgn}(D_0))$
is an odd $p$-summable bounded Breuer-Fredholm
module. Now, each $p$-summable
bounded Breuer-Fredholm module for $\cal A$
has associated with it a cyclic $(p-1)$-dimensional cycle
over $\cal A$
(see [Co2] p.292). We may  therefore utilise the standard formula for the
character of this
 cycle which is, with $p=2n+2$,
$$\tau_{2n+1}(a^0,a^1,\ldots,a^{2n+1})
=\tau(F_0[F_0,a^0][F_0,a^1]\ldots[F_0,a^{2n+1}])$$
where $a^j\in \cal A$. Thus $\tau_{2n+1}$
may be regarded as
the cyclic cocycle associated with both the odd unbounded
Breuer-Fredholm module and the bounded one.

Note that in [Co1] in the case of even $p$-summable
Fredholm modules an expression
is given directly in terms of the unbounded operator $D_0$ for an
associated cyclic cocycle. We have not investigated the existence
of such an expression in the odd case but presumably one exists.

It is worth mentioning at this point the
 motivating example of spectral flow in [CP].
Given an odd $p$-summable unbounded Breuer-Fredholm module
$({\cal M}, D_0)$ [CP] (following [G]) introduce an affine space
$\Phi_p = \{D=D_0+A\ | \ A\in {\cal M}_{sa}\}$. Then it
is shown in [CP] that the map
$D\mapsto \phi(D)= D({\bf 1}+D^2)^{-1/2}$ takes, for all $q>p$,
the space $\Phi_p$ continuously into the affine space
$${\cal M}_q= \{F=F_0+ X | \ X\in  L_{q,q/2}({\cal M},\tau)_{sa}\}$$
where $F_0=\phi(D_0)$ and
 $ L_{q,q/2}({\cal M},\tau)_{sa}$ are the bounded
self adjoint elements
of  $ L_q({\cal M},\tau)$ that satisfy the additional condition
$$XF_0+F_0X\in L_{q/2}({\cal M},\tau).$$
This second constraint plays a key role in the
analytic formulae for spectral flow in [CP].
In fact the import of Corollary 0.4 is that $\Phi_p$
actually maps continuously into ${\cal M}_p$
by the following Proposition.

\noindent {\bf  Proposition 3.2 \quad}
{\it If $1<p<\infty$ and $({\cal M}, D_0)$ is an odd unbounded
$p$-summable Breuer-Fredholm module then the map
$$\phi: \Phi_p \rightarrow \phi(D_0) + L_{p,p/2}({\cal M},\tau)_{sa}$$
is well defined and continuous.}

\noindent{\bf Proof. \quad}
By Theorem 0.3(i) $\phi(D)=F$ lies in  
$\phi(D_0)+ L_{p}({\cal M},\tau)_{sa}$
and the mapping is continuous in that space.
To take account of the additional condition
set $F_0=\phi(D_0)$ and define $X_D=F-F_0$ so that the map
$D\mapsto X_D\in L_{p}({\cal M},\tau)_{sa}$
is continuous on $\Phi_p$. 
Now the map 
$$D\mapsto 1-\phi(D)^2=(1+D^2)^{-1}\in L_{p/2}({\cal M},\tau)_{sa}$$
is continuous by Corollary A.2 (or [CP], Proposition 10 of Appendix B).
But 
$$D\mapsto X_D\mapsto X_D^2\in L_{p/2}({\cal M},\tau)_{sa}$$
is also continuous. We observe that 
$$1-\phi(D)^2=1-F_0^2 -(X_D^2 +F_0X_D+X_DF_0)= (1+D_0^2)^{-1}
-X_D^2-(F_0X_D+X_DF_0)$$
and so $D\mapsto(F_0X_D+X_DF_0)\in  L_{p/2}({\cal M},\tau)_{sa}$
is continuous. That is $X_D \in L_{p,p/2}({\cal M},\tau)_{sa}$
and continuity in the norm $||.||_{p,p/2}$ on the latter space is clear from
the definition of this norm:
$$||X_{D_1}-X_{D_2}||_{p,p/2}
=||X_{D_1}-X_{D_2}||_p +|| (X_{D_1}-X_{D_2})F_0 +F_0(X_{D_1}-X_{D_2})||_{p/2}.
\quad \square$$

In particular continuous paths
$$\{D_t= D_0+A_t\}\eqno (3.3)$$
in $\Phi_p$ map to continuous paths of Breuer-Fredholm operators
in ${\cal M}_p$ under $\phi$. Thus using [P2] we
define the spectral flow along the path $\{D_t\}$
as the spectral flow along $\{\phi(D_t)\}$ and this will be
independent of the path in ${\cal M}_p$
joining the endpoints as ${\cal M}_p$ is simply connected.

\noindent {\bf Remarks 3.3\quad }
Generalizing  the notions of
$p$-summable and $\theta$-summable Fredholm modules (see [Co1,2])
one of us  introduced
in [S] the notion of an odd (un)bounded Breuer-Fredholm module 
associated with
an arbitrary symmetric operator space $\em$ (which we now abbreviate to
$\em$-summable Breuer-Fredholm module). For the definitions
and additional information concerning these spaces we refer to [DDP],
[DDPS], [SC]. An inspection of the proofs presented in Section 2
shows that the assertions given in Proposition 2.1, 2.2, 2.3 and 
Theorem 0.3
would also hold if

\noindent (i) $\lp$ is replaced by an arbitrary symmetric 
operator space $\em$ which is an interpolation space for any couple 
$(L_{p_{_{1}}}(\nm, \tau ), L_{p_{_{2}}}(\nm ,\tau ))$ with
$1<p_1\leq p_2<\infty $ and 

\noindent(ii) which  has the Fatou property (see e.g.
[DDPS]). 

\noindent In particular, these conditions are satisfied when the corresponding
symmetric function space $E$ has non-trivial Boyd indices 
(see e.g. [LT]) and the Fatou property (see [DDP], [DDPS]). 
The latter assumption about the Fatou property is automatically 
satisfied whenever $E$ is an Orlicz, Lorentz or Marcinkiewicz function
space (see [LT], [BS], [KPS]).

\noindent By way of an example, it
 follows from  [FG], Theorem 4.1 that any reflexive Orlicz space 
$L_\Phi$ has non-trivial Boyd indices. It is also
well-known that the latter property holds also for the family of
spaces $L_{_{p,q}},\ 1<p<\infty,\ 1\leq q\leq \infty$ (for the
definition of the latter spaces we refer to [LT]; the spaces
$L_{_{p,\infty }}$ are known as {\it weak $L_p$-spaces} see e.g.
[BS]). Thus we obtain the following strengthening of Corollary 0.4.

\noindent {\bf Corollary ${\bf 0.5}$\quad } {\sl If $\em$ is either 
a reflexive Orlicz operator space $L_\Phi(\nm ,\tau )$, or  
 $L_{_{p,q}}(\nm ,\tau ),\ 1<p<\infty,\ 1\leq q\leq \infty$
and if $(\nm , D_0)$ 
is an odd unbounded $\em$-summable 
Breuer-Fredholm module for ${\cal A}$,
 then $(\nm , {\rm sgn}(D_0))$ is an odd  bounded  $\em$-summable 
Breuer-Fredholm module for
${\cal A}$.}

\noindent It is also worth mentioning specifically
 that this last Corollary holds when $\em$
is the space  $L_{_{p,\infty}}(\nm ,\tau ),$ $1<p<\infty$.
This 
latter fact may also be deduced from [CP] and [CM].
It has implications for analytic formulae for spectral
flow.

\bigskip
\noindent{\bf Appendix}.

\noindent
The result we need in theorem 0.3(i) is the following.

\noindent {\bf Theorem \quad } {\sl  Let
 $({\cal M}, \tau)$ be an arbitrary semifinite
von Neumann algebra,
let $x=x^*$ be affiliated with ${\cal M}$ such
 that $z=({\bf 1} + x^2)^{-1/2}\in \lp$. Then
there exists a constant $C>0$ (depending on $\lp$ and $x$) such that  for
all bounded self-adjoint $y-x\in {\cal M}$ we have }
$$
\biggl \Vert {|y|\over (1+y^2)^{1/2}}-{|x|\over (1+x^2)^{1/2}}
\biggr \Vert _{_{\lp}} \leq C\max \{\Vert y-x\Vert^{1/2},\Vert y-x\Vert\}.
\eqno (A.1)
$$

\bigskip
We introduce some notation.
We let $\mu(x)$ denote the generalised singular value function
for $x\in\nmt$ (see [FK] for details).
 If $x,y\in \nmt$, then we say that
$x$ is submajorized by $y$ and write $x\prec\prec y$ if and only if
$$
\int _{_{0}}^t\mu _{_{s}}(x)ds
\leq \int _{_{0}}^t\mu _{_{s}}(y)ds,\quad t\ge 0.
$$

\noindent The proof of the theorem rests on the following:

\bigskip
\noindent {\bf Lemma A.1\quad } {\sl Let $x=x^*\in \nmt$ and
 $0\leq y \in \nmt$
and let $-y\leq x\leq y$. Then $\mu_s(x)\leq 
\mu_{s/2}(y)$ for all $s>0$}.

\bigskip
\noindent {\bf Proof \quad } We let $p_\pm$ denote the 
spectral projections corresponding to
  the positive and negative parts of the spectrum of $x$.
We have $x_+\leq p_+yp_+$ and $x_-\leq  p_-yp_-$.
So $\mu_s(x_+)\leq \mu_s(y)$ and $\mu_s(x_-)\leq \mu_s(y)$, whence
$\mu_s(x_+^n)\leq \mu_s(y)^n$ and $\mu_s(x_-^n)\leq \mu_s(y)^n$ 
for all $n=1,2,\dots$ and all $s>0$. Thus
$$
\mu_s(x)^n=\mu_s(|x|^n)=\mu_s(x_+^n+x_-^n)\leq\mu_{s/2}(x_+^n)+
\mu_{s/2}(x_-^n)\leq 2 \mu_{s/2}(y)^n, \quad n=1,2, \dots \ .
$$
It follows that $\mu_s(x)\leq \mu_{s/2}(y),\ s>0. \quad \square $

\bigskip
\noindent{\bf Corollary A.1.} {\sl If $f$ is any
continuous increasing function
on $[0,\infty)$ with $f(0)\geq 0$ and  $-y\leq x\leq y$ then
$\mu_s(f(|x|)\leq \mu_{s/2}(f(y))$
and
$$\int_0^r \mu_s(f(|x|))ds \leq 2 \int_0^r\mu_s(f(y))ds.$$
That is, $f(|x|)\prec\prec 2f(y)$ so that in particular
$|x|^{1/2}\prec\prec 2y^{1/2}$.}

\bigskip
\noindent {\bf Proof \quad } The first claim follows from [FK] lemma 2.5(iv).
The second claim holds because:
$$\int_0^r \mu_s(f(|x|))ds \leq\int_0^r\mu_{s/2}(f(y))ds=
2\int_0^{r/2}\mu_t(f(y))dt\leq 2 \int_0^r\mu_t(f(y))dt.\quad \square$$

\bigskip

\noindent{\bf Corollary A.2.} 
{\sl Given a fully symmetric function space $\em$ and
 any continuous increasing function
$f$ on $[0,\infty)$ with $f(0)\geq 0$ we have,
for $x=x^*$,$y\geq 0$, $-y\leq x\leq y$ and 
$f(y)\in E({\cal M},\tau)$,
that $f(|x|)\in \em$ and
$||f(|x|)||_{E({\cal M},\tau)}\leq 2||f(y)||_{E({\cal M},\tau)}$. }

\bigskip 

\noindent We now move to the proof of the theorem.

\bigskip 
\noindent {\bf Proof \quad }  It follows from the proof in [CP]
Appendix B, Proposition 10 that
$$
-2\max \{\Vert y-x\Vert^{2},\Vert y-x\Vert\} \cdot
{1\over (1+x^2)}\leq
{1\over (1+y^2)}-{1\over (1+x^2)}\leq
2\max \{\Vert y-x\Vert^{2},\Vert y-x\Vert\} \cdot {1\over (1+x^2)}
$$
and, by Corollary A.1, it follows that
$$
\biggl |{1\over (1+y^2)}-{1\over (1+x^2)}\biggr |^{1/2} \prec\prec
2( 2\max \{\Vert y-x\Vert^{2},\Vert y-x\Vert\} )^{1/2}\biggl ({1\over 
(1+x^2)^{1/2}}\biggr )\eqno (A.2).
$$
The latter inequality implies immediately that
$$
\biggl \Vert\ \biggl |{1\over (1+y^2)}-{1\over (1+x^2)}\biggr |^{1/2}
\biggr \Vert _{_{\lp}}
\leq  2^{3/2}\max \{\Vert y-x\Vert^{1/2},\Vert y-x\Vert\} \biggl \Vert
 {1\over (1+x^2)^{1/2}}
\biggr\Vert _{_{\lp}}
$$
or, equivalently,
$$
\biggl \Vert\ \biggl |{y^2\over (1+y^2)}-{x^2\over (1+x^2)}\biggr |^{1/
2}\biggr \Vert _{_{\lp}} \leq  C\max \{\Vert y-x\Vert^{1/2},\Vert y-x\Vert\}
  \eqno (A.3)
$$
where $C:=2^{3/2}\biggl \Vert {1\over (1+x^2)^{1/2}} \biggr
\Vert _{_{\lp}}$.
\smallskip
\noindent
The following inequality is developed in [BKS] for the case of
symmetrically-normed ideals of compact operators,
extended to measurable operators affiliated with an
arbitrary semifinite von Neumann algebra ${\cal M}$ by H. Kosaki
(it is given in the appendix to [HN]) with an alternative version of the
proof given in [DD]. By Theorem 1.1 from [DD]
(see also [BKS] Theorem 1) we have
$$
a^{1/2}-b^{1/2} \prec \prec |a-b|^{1/2} \eqno (A.4)
$$
for any $0\leq a,b\in \nmt$. Combining (A.4) with (A.3) and
using the fact that  the space
$\lp$ is fully symmetric we have
$$
\eqalign {
\biggl \Vert {|y|\over (1+y^2)^{1/2}}-
{|x|\over (1+x^2)^{1/2}}\biggr \Vert _{_{\lp}}
&=\biggl \Vert \biggl ({y^2\over (1+y^2)}\biggr )^{1/2}-
\biggl ({x^2\over (1+x^2)}\biggr )^{1/2}\biggr \Vert _{_{\lp}} \cr
&\leq \biggl \Vert\ \biggl |{y^2\over (1+y^2)}-{x^2\over (1+x^2)}\biggr
 |^{1/2}\biggr \Vert _{_{\lp}} \cr
&\leq C\max \{\Vert y-x\Vert^{1/2},\Vert y-x\Vert\}\cr }
$$
which is the desired inequality.
\quad$\square$

\bigskip

\bigskip
\bigskip

\parindent=40pt
\parskip=6pt plus1pt
\def\bib#1{\item{\hbox to \parindent{#1\hfill}}}
\overfullrule=0pt
\baselineskip=12pt

\bigskip
\bigskip

\centerline {R E F E R E N C E S}

\bigskip



\bib {[BKS]} {M.S. Birman, L.S. Koplienko and M.Z. Solomyak, {\it Estimates
 for the spectrum of the difference between fractional powers of two
 selfadjoint operators,} Soviet Mathematics, {\bf 19}(3) (1975), 1-6.}



\bib {[BS]} {C. Bennett and R. Sharpley, {\it Interpolation of operators,}
 Pure and applied mathematics, v.129 (1988).}

\bib {[BGM1]} {E. Berkson, T.A. Gillespie and P.S. Muhly, {\it Abstract
spectral decomposition guaranteed by the Hilbert transform,}
Proc. London Math. Soc. {\bf 53} (1986), 489-517.}

\bib {[BGM2]} {E. Berkson, T.A. Gillespie and P.S. Muhly, {\it
Generalized analyticity in $UMD$-spaces,} Arkiv f\"or Math {\bf
27}(1989), 1-14.}

\bib {[BSo1-3]} {M.Sh. Birman and M.Z. Solomyak, {\it Double Stieltjes operator
integrals $I;II;III$}, Problemy Mat. Fiz. {\bf 1} (1966), 33-67; {\bf 2}
(1967), 26-60; {\bf 6} (1973), 27-53.}

\bib {[CP]}  { A. Carey and J. Phillips, {\it Unbounded Fredholm modules
and spectral flow,} Canadian J. Math. {\bf 50} (1998), 673-718.}

\bib {[CP1]} {A. Carey and J. Phillips, preprint in preparation} 

\bib {[CM]} {A. Connes and H. Moscovici, {\it
Transgression du caractere de Chern et cohomologie cyclique},
C.R. Acad. Sc. Paris Ser. I Math. {\bf 303} (1986)
913-918.}

\bib {[Co1]} {A. Connes, {\it Noncommutative differential geometry},  Publ.
Math. Inst. Hautes Etudes Sci. {\bf 62} (1985), 41-144.}

\bib {[Co2]} { A. Connes, {\it Noncommutative geometry}, Academic Press, (1994).}


\bib {[Da1]} {E.B. Davies, {\it One parameter semigroups}, Academic Press,
 London (1980).}

\bib {[Da2]} {E.B. Davies, {\it Lipschitz continuity of
functions of operators in the Schatten classes}, J. London Math. Soc.
{\bf 37} (1988), 148-157.}

\bib {[D]} {J. Dixmier, {\it Von Neumann algebras}
Elsevier, Amsterdam, 1981.}

\bib {[DD]}  { P.G. Dodds, T.K. Dodds, {\it On a submajorization inequality of
T. Ando,} Operator Theory Advances and Applications {\bf 75} (1995), 113-133.}



\bib {[DDP]} {P.G. Dodds, T.K. Dodds and B. de Pagter, {\it
Non-commutative  K\"othe duality}, Trans. Amer. Math. Soc. {\bf 339}
(1993), 717-750.}

\bib {[DDPS]}  { P.G. Dodds, T.K. Dodds, B. de Pagter and F.A. Sukochev,
{\it Lipschitz continuity of the absolute value and Riesz projection in
symmetric operator spaces}, J. Functional Anal. {\bf 148} (1997), 28-69.}




\bib  {[FG]} {D.L.Fernandez and J.B.Garcia, {\it Interpolation of
Orlicz-valued function spaces and UMD property}, Studia Math. {\bf 99}
(1991), 23-40.}

\bib {[FK]} {T. Fack and H. Kosaki, {\it Generalized s-numbers of
$ \tau $-measurable operators}, Pacific J. Math.
{\bf 123} (1986), 269-300.}

\bib {[G]} {E. Getzler, {\it The odd Chern character in cyclic homology 
and spectral flow}, Topology {\bf 32} (1993), 489-507.}

\bib {[HN]} {F. Hiai and Y. Nakamura, {\it Distance between unitary orbits in
 von Neumann algebras,} Pacific J. Math. {\bf 138} (1989), 259-294.}


\bib {[KPS]} {S.G. Krein, Ju.I. Petunin and E.M.
Semenov, {\it Interpolation of linear operators}, Translations of
Mathematical Monographs, Amer. Math. Soc. {\bf 54} (1982).}


\bib {[LT]} {J. Lindenstrauss and L. Tzafriri, {\it
 Classical Banach Spaces II}, Springer-Verlag, (1979).}

\bib {[Ne]} {E. Nelson, {\it Notes
on non-commutative integration}, J. Functional Anal. {\bf 15} (1974),
103-116.}



\bib {[Pa]} {T.W. Palmer, {\it $*$-representations of $U^*$-algebras},  
Indiana Univ. Math. J. {\bf 20B} (1971), 929-933.}

\bib {[P1]} {J. Phillips, {\it Self-adjoint Fredholm operators and spectral flow,} Canad. Math. Bull. {\bf 39} (1996), 460-467.}

\bib {[P2]} {J. Phillips, {\it Spectral flow in type $I$ and $II$ factors - A new approach}, Fields Inst. Comm. {\it Cyclic Cohomology $\&$ Noncommutative Geometry}, {\bf 17} (1997), 137-153.}

\bib {[RS1]} {M. Reed and B. Simon, {\it Methods of modern mathematical
physics I},  Academic Press, New York and London (1972).}

\bib {[RS2]} {M. Reed and B. Simon, {\it Methods of modern mathematical
physics II},  Academic Press, New York and London (1972).}

\bib {[RN]} {F. Riesz and B. Sz.-Nagy {\it Functional Analysis} Frederick Ungar
Publishing Co., 1955.}

\bib  {[Si]} { B. Simon, {\it Trace ideals and their applications}, 
London Math. Soc. Lecture Notes Series {\bf 35} (1979).}

\bib  {[S]} {F.A. Sukochev , {\it Operator esimates for Fredholm
modules}, preprint.}

\bib  {[SC]} {F.A. Sukochev and V.I. Chilin, {\it Symmetric spaces on
semifinite von Neumann algebras,} Soviet Math. Dokl {\bf 42} (1991), 97-101.}

\bib  {[SZ]} { S. Stratila and L. Zsido, {\it Lectures on von Neumann algebras}, Abacus Press (1979).}

\bib {[Ta]} {M. Takesaki,\  {\it Theory of Operator Algebras I},\
Springer-Verlag,\  New York-Heidelberg \  Berlin, 1979.}

\bib {[Te]} {M. Terp, {\it $L^p$-spaces associated with von Neumann
algebras,} Notes, Copenhagen Univ. (1981).}

\bib {[Ti]} {O.E. Tikhonov, {\it Continuity of operator functions in
topologies connected to a trace in  Neumann's algebra} (Russian),
Izv. Vyssh. Uchebn. Zaved. Mat., 1987, no.1, 77-79; translated in Soviet
Math. (Iz.VUZ), {\bf 31} (1987), 110-114.}


\end

which is based on the use of reflexivity of the space
$\lp,\ p\in (1,\infty)$.
It follows immediately from (2.12) that for an arbitrary projection
$q\in \nm$ we have
$$
q({e_nye_n\over ({\bf 1}+(e_nye_n)^2)^{1/2}} -
{e_nxe_n\over ({\bf 1}+(e_nxe_n)^2)^{1/2}} )q
\ \ {\buildrel (so)\over \longrightarrow}\ \
q({y\over ({\bf 1}+y^2)^{1/2}} -
 {x\over ({\bf 1}+x^2)^{1/2}})q,\quad  n\to \infty. \eqno (2.17)
$$
It follows from (2.17) that
$$
\tau (q({e_nye_n\over ({\bf 1}+(e_nye_n)^2)^{1/2}} -
{e_nxe_n\over ({\bf 1}+(e_nxe_n)^2)^{1/2}} )q)\to
\tau (q({y\over ({\bf 1}+y^2)^{1/2}} -
 {x\over ({\bf 1}+x^2)^{1/2}})q). \eqno (2.18)
$$
since the restriction of the trace to the unit ball
of the algebra $q{\cal M}q$ is strong operator continuous (see [SZ] p.16).
Furthermore the sequence
$$
\{{e_nye_n\over ({\bf 1}+(e_nye_n)^2)^{1/2}} -
{e_nxe_n\over ({\bf 1}+(e_nxe_n)^2)^{1/2}}\}_{_{n=1}}^\infty \subseteq \lp
$$
is norm bounded (see (2.10)) so that reflexivity of $\lp$
implies  there exists an element
$w\in \lp$ and a subsequence $(n_k)$ such that
$$
\sigma (\lp, L_{p'}(\nm ,\tau ))-\lim _{k\to \infty }
({e_{n_{_{k}}}ye_{n_{_{k}}}\over
({\bf 1}+(e_{n_{_{k}}}ye_{n_{_{k}}})^2)^{1/2}} -
{e_{n_{_{k}}}xe_{n_{_{k}}}\over
({\bf 1}+(e_{n_{_{k}}}xe_{n_{_{k}}})^2)^{1/2}}) =w,\eqno (2.19)
$$
where $L_{p'}(\nm ,\tau )=\lp^*$. In particular, it follows from (2.19) that
$$
\tau (q({e_{n_{_{k}}}ye_{n_{_{k}}}\over
({\bf 1}+(e_{n_{_{k}}}ye_{n_{_{k}}})^2)^{1/2}} -
{e_{n_{_{k}}}xe_{n_{_{k}}}\over
({\bf 1}+(e_{n_{_{k}}}xe_{n_{_{k}}})^2)^{1/2}} )q)\to
\tau (qwq),\eqno (2.20)
$$
as $k\to \infty$. Combining (2.18) and (2.20) we conclude that
$$
\tau (q({y\over ({\bf 1}+y^2)^{1/2}} -
 {x\over ({\bf 1}+x^2)^{1/2}})q)=\tau (qwq), \quad {\rm or}\quad
\tau (({y\over ({\bf 1}+y^2)^{1/2}} -
 {x\over ({\bf 1}+x^2)^{1/2}})q)=\tau (wq)
$$
for any projection $q\in \nm$ with $\tau (q)<\infty$. It follows immediately
that
$$
w={y\over ({\bf 1}+y^2)^{1/2}} -
 {x\over ({\bf 1}+x^2)^{1/2}}. \eqno (2.21)
$$
It follows from (2.19) that
$$
\Vert w\vlp\leq \liminf_{_{n}} \Vert {e_{n_{_{k}}}ye_{n_{_{k}}}\over
({\bf 1}+(e_{n_{_{k}}}ye_{n_{_{k}}})^2)^{1/2}} -
{e_{n_{_{k}}}xe_{n_{_{k}}}\over
({\bf 1}+(e_{n_{_{k}}}xe_{n_{_{k}}})^2)^{1/2}}\vlp
$$
and combining now (2.10) and (2.21) we conclude that
$$
\Vert {y\over ({\bf 1}+y^2)^{1/2}} -
 {x\over ({\bf 1}+x^2)^{1/2}}
\vlp
\leq c_p'\max \{ \Vert x-y\Vert ^{1/2}, \ \Vert x-y\Vert \}.
$$
Letting  (see equality (2.11))
$$
{\cal Z}_p':={\cal K}_p+2^{3/2}+1
$$
we arrive to the inequality $(0.3)'$.
The inequality (0.3) of Theorem 0.3(i) follows from the inequality $(0.3)'$
 via Proposition 2.3.\quad $\square$

In the norm
$||a||_0 = K||a||_{\cal A} + ||[D_0,a]||$
the constant $K \geq 1$  is included to ensure that
the norm is submultiplicative, that is $||ab||_0\leq||a||_0||b||_0$. To
see that ${\cal A}_0$ is invariant under $*$,
let $a\in {\cal A}_0$
and $\xi\in$ dom$D_0$ then $a^*\xi\in$ dom$D_0$
if and only if $\eta\rightarrow <a^*\xi,D_0\eta>$
is bounded for $\eta\in$ dom$D_0$.
But
$$<a^*\xi,D_0\eta>=<\xi,aD_0\eta>=<\xi, ([a,D_0]+ D_0a)\eta>
=<\xi, [a,D_0]\eta> +<D_0\xi, a\eta>$$
which gives a bounded linear functional for  $\eta\in$ dom$D_0$.
One then checks that $[D_0,a^*]=-[D_0,a]^*$ so that $a^*\in{\cal A}_0$.
To see that ${\cal A}_0$ is complete in $||.||_0$
norm consider a Cauchy sequence $\{a_n\}$. 
Then $a_n\rightarrow a$ and $[D_0,a_n]\rightarrow b$ where
$a$ and $b$ are
bounded operators.
Using the same idea as above we consider for $\xi,\eta\in$ dom$D_0$
$$<a\xi,D_0\eta>=\lim_{n\rightarrow\infty}<\xi,a_n^*D_0\eta>
=\lim_{n\rightarrow\infty}<\xi, ([a^*_n,D_0]+ D_0a^*_n)\eta>$$
$$
=\lim_{n\rightarrow\infty}(<\xi, -[D_0,a_n^*]\eta>+<D_0\xi,a^*_n\eta>)
=<\xi, -b^*\eta>+<D_0\xi,a^*\eta>$$
Thus we conclude that 
$a$ preserves the domain of $D_0$.
{}From this is is now easy to see that $b= [D_0,a]$
and that the sequence $\{a_n\}$
has limit $a$ in $||.||_0$-norm.$\quad\square$

Thus the import of Corollary 0.4 is that $\Phi_p$
actually maps continuously into ${\cal M}_p$.

We now sketch the implications for
 spectral flow formulae (full details depend on
[CP] and work in progress [CP1]).
In [CP] spectral flow was given in terms of the integral of
a one form on the tangent space to $\Phi_q$
for integers $q>p$. For non-integer
$q>p$ we refer to [CP1]. Thus one writes
for $D\in\Phi_q$,
$$\alpha_q(A)={{1}\over{\tilde C_{q/2}}}\tau(A({\bf 1}+D^2)^{-q/2})$$
where $A\in {\cal M}_{sa}$ is a tangent
vector at $D$ and $\tilde C_{q/2}$ is a normalisation constant.
An analogous one form was introduced on the tangent space to
${\cal M}_q$ for all $q>p$
by setting:
$$\theta_q(X)= {{1}\over{C_{(q-1)/2}}}
\tau(X({\bf 1}-\phi(D_0)^2)^{{q-1}\over{2}})$$
for $X\in  L^{q,q/2}({\cal M},\tau)_{sa}$ and $C_{(q-1)/2}$
a normalisation constant.
Now the one form
$$\theta_p(X)= {{1}\over{C_{(p-1)/2}}}
\tau(X({\bf 1}-\phi(D_0)^2)^{{p-1}\over{2}}) $$
is also well defined for  $X\in  L^{p,p/2}({\cal M},\tau)_{sa}$.
The one forms $\theta_q$ for $q>p$
are used in [CP] and [CP1] to give analytic
expressions for spectral flow along piecewise differentiable
paths in $\Phi_p$ and ${\cal M}_q$ for $q>p$.
(Note that the prototype for such analytic expressions appeared
in [G].)
More precisely, it follows from [CP] that for all $q\geq p$
(provided $(q-1)/2$ is an integer)
and all continuous
paths $\{F_t\}$
of bounded self-adjoint Breuer-Fredholm operators in ${\cal M}_q$
which are piecewise $C^1$,
and whose endpoints are unitarily equivalent,
the spectral flow is given
by
$$sf\{F_t\}={{1}\over{C_q}}\int_1^2 \tau
    ({{d}\over{dt}}(F_t)({\bf 1}-F_t^2)^{{q-1}\over{2}}) dt.  \eqno (3.4)$$
When $F_t= \phi(D_t)$ with $D_t$ as in (3.3) 
with unitarily equivalent endpoints then (3.4) 
follows for $q>p$ (and $(q-1)/2$ an integer) by [CP].
In order to establish (3.4) for $q=p$
(whether or not  $(q-1)/2$ an integer) we argue as follows
provided $\{ \phi(D_t)\}$ is a piecewise $C^1$
path with tangent vectors in
$L^{p,p/2}({\cal M},\tau)_{sa}$.
First, $\phi(D_t)\in{\cal M}_p$  by (A.2) of the appendix.
Thus the integrand in (3.4) is well defined for $q=p$.
We then generalise [CP] so that (3.4) holds
for $q>p$ even when $(q-1)/2$ is not an integer
(see the forthcoming [CP1]). Finally we
take the limit as $q\to p$ in the integral formula
to obtain the result.

\noindent Note that there is an alternative unbounded
spectral flow formula established in [CP] and [CP1] using the forms
$\alpha_q$:
$$sf\{D_t\}={{1}\over{\tilde C_m}}\int_1^2 \tau ({{d}\over{dt}}(D_t)
(({\bf 1}+D_t^2)^{-m}) dt$$
for $m>1+p/2$ and a certain constant $\tilde C_m$.
The relationship between this and (3.3) is given by the differentiation 
formula
$$\tau({{d}\over{dt}} \phi(D_t)({\bf 1}- \phi(D_t)^2)^{n})
=\tau({{d}\over{dt}}(D_t)({\bf 1}+D_t^2)^{-n-3/2})$$
for $n> (p-1)/2$. Theorem 0.3(i) is not sufficient to 
establish this formula for piecewise differentiable paths
$\{D_t\}$ when $n=(p-1)/2$.
 Thus the exponents in the bounded and unbounded
spectral flow formulae are not optimal. Sharper results are expected from
work in progress using ideas from the study of
spectral flow in  $\theta$-summable Breuer-Fredholm modules.